\documentclass[11pt]{article}          
\usepackage{latexsym}
\usepackage{amssymb}
\usepackage{bbm}
\usepackage[all]{xy}
\usepackage{graphicx}
\usepackage{amsmath}

\makeatletter
\def\@begintheorem#1#2{\trivlist%
   \item[\hskip \labelsep{\sffamily\bfseries #1\ #2}]\itshape}
\newtheorem{teo}{Theorem}[section]
\newtheorem{defi}[teo]{Definition}

\newtheorem{pro}[teo]{Proposition}
\newtheorem{_rem}[teo]{Remark}
\newtheorem{_eje}[teo]{Example}
\newtheorem{_conj}[teo]{Conjecture}
\newenvironment{rem}{\def\@begintheorem##1##2{\trivlist%
 \item[\hskip\labelsep{\sffamily\bfseries ##1\ ##2}]}\begin{_rem}}{\end{_rem}}

\makeatother

\newenvironment{beweis}{{\em Proof:}}{\hfill $\rule{2mm}{2mm}$
\vspace{3mm}

}

\DeclareMathAlphabet{\Ma}{U}{msa}{m}{n}
\DeclareMathAlphabet{\Mb}{U}{msb}{m}{n}
\DeclareMathAlphabet{\Meuf}{U}{euf}{m}{n}
\DeclareSymbolFont{ASMa}{U}{msa}{m}{n}
\DeclareSymbolFont{ASMb}{U}{msb}{m}{n}
\DeclareMathSymbol{\hrist}{\mathord}{ASMa}{"16}
\DeclareMathSymbol{\varkappa}{\mathalpha}{ASMb}{"7B}
\DeclareMathSymbol{\CrPr}{\mathord}{ASMb}{"6F}

\def\got#1{\Meuf{#1}}

\def\mr #1.{\mathrm{#1\,}}
\def\mrt #1.{\mathrm{\mbox{\tiny #1\,}}}   
\def\mt #1.{{\mbox{\tiny $#1$}}}
\def\ms #1.{{\mbox{\small $#1$}}}
\def\ol #1.{\overline{#1}}
\def\mb #1.{\mathbf{#1\,}}
 
\def\1{\mathbbm 1}

  \def\al #1.{{\mathcal{#1}}}
  \def\ot #1.{{\got{#1}}}
  \def\C{\Mb{C}}
  
  \def\kl #1.{{\mbox{\tiny {$\Mb{#1}$}}}}

  \def\N{\Mb{N}}
  \def\Z{\Mb{Z}}

  \def\t #1.{\tilde{#1}} 
  \def\T #1.{\widetilde{#1}}

\DeclareMathSymbol{\hsemi}{\mathord}{ASMb}{"6E}
\newcommand{\semi}[2]{\mbox{$#1\kern.1em\hsemi\kern.1em#2$}}

\def\LA{\left\langle\bgroup}
\def\LE{\left[\bgroup}
\def\LG{\left\{\bgroup}
\def\LR{\left(\bgroup}
\def\RA{\egroup^{\rule{0mm}{2mm}}\right\rangle}
\def\RE{\egroup^{\rule{0mm}{2mm}}\right]}
\def\RG{\egroup^{\rule{0mm}{2mm}}\right\}}
\def\RR{\egroup^{\rule{0mm}{2mm}}\right)}
\def\Ldummy{\left.\bgroup}
\def\Rdummy{\egroup^{\rule{0mm}{2mm}}\right.}

\def\Kbegin{\begin{equation} \left. \begin{array}{rcl}}
\def\Kend{\end{array} \right\} \end{equation}}

\def\l2{\Lambda^{\mbox{\tiny $(2)$}}}

  \def\ccr #1,#2.{\overline{\Delta(#1,\,#2)}}

  \def\b #1.{{\bf #1}}
  \def\cross#1.{\mathrel{\mathop{\times}\limits_{#1}}}
  \def\C{\Mb{C}}
  \def\N{\Mb{N}}

  \def\Z{\Mb{Z}}
  \def\T{\Mb{T}}

  \def\wwh #1.{\widehat{#1}}
  \def\wt #1.{\widetilde{#1}}

  \def\cross #1.{\mathrel{\raise 3pt\hbox{$\mathop\times\limits_{#1}$}}}
\def\set #1,#2.{\left\{\,#1\;\bigm|\;#2\,\right\}}
\def\b #1.{{\bf #1}}

\def\aut{{\rm Aut}\,}

\def\ol #1.{\overline{#1}}
\def\rn#1.{\romannumeral{#1}}

\def\s #1.{_{\smash{\lower2pt\hbox{\mathsurround=0pt $\scriptstyle #1$}}\mathsurround=3pt}}
\def\bra #1,#2.{{\left\langle #1,\,#2\right\rangle_{\al A.}}}

\def\XP#1!{\renewcommand{\baselinestretch}{.7}\marginpar{{\footnotesize #1}\hfil}
\renewcommand{\baselinestretch}{1.5}}
\def\XB{\marginpar{
{\footnotesize\bf Change~starts----}\lower 11pt\hbox{\mathsurround=0pt$
\!\!\displaystyle{
\Bigg\downarrow}$\mathsurround=3pt}}}
\def\XE{\marginpar{{\footnotesize\bf Change~ends-----}\raise 10pt
\hbox{\mathsurround=0pt$ 
\!\!\displaystyle{
\Bigg\downarrow}$\mathsurround=3pt}}}

\title{\bf On the nuclearity of certain 
           Cuntz-Pimsner algebras}
\author{           
 {\sc Fernando Lled\'o}  \\[2mm] 
 {\footnotesize Institute for Pure and Applied Mathematics,}       \\
 {\footnotesize RWTH-Aachen, Templergraben 55,}                   \\ 
 {\footnotesize D-52062 Aachen, Germany.}                         \\[1mm]
 {\footnotesize lledo@iram.rwth-aachen.de}
\and
 {\sc Ezio Vasselli}\\[2mm] 
 {\footnotesize Dipartimento di Matematica,}\\
 {\footnotesize University of Rome "La Sapienza"}  \\       
 {\footnotesize P.le Aldo Moro 2, I-00185 Roma, Italy}    \\[1mm]
 {\footnotesize vasselli@mat.uniroma2.it}}

\date{\today{}}
\begin{document}
\maketitle

\begin{abstract} 
In the present paper, we give a short proof of the nuclearity 
property of a class of Cuntz-Pimsner algebras associated 
with a Hilbert $\al A.$-bimodule $\al M.$, 
where $\al A.$ is a separable and nuclear {\it C*}-algebra.
We assume that the left $\al A.$-action on the bimodule $\al M.$
is given in terms of compact module operators and that $\al M.$ is 
direct summand of the standard Hilbert module over $\al A.$.
\end{abstract}


\section{Introduction}

M.V.~Pimsner introduced in his seminal paper
\cite{Pimsner97} a new family of {\it C*}-algebras 
$\al O._\al M.$ that are naturally generated by a Hilbert bimodule
$\al M.$ over a {\it C*}-algebra $\al A.$. These algebras
generalise Cuntz-Krieger algebras as well as crossed-products by the group
$\Z$. In Pimsner's construction $\al O._\al M.$ is given as
a quotient of a Toeplitz like algebra acting on a concrete
Fock space associated to $\al M.$. 
An alternative abstract approach to Cuntz-Pimsner algebras
in terms of {\it C*}-categories is given in 
\cite{Doplicher98,Kajiwara98,PinzariIn97}
(for the notion of C*-category see \cite{Doplicher89b}).

In the present note, we give a short proof of nuclearity 
property of the Cuntz-Pimsner
algebra (cf.~Theorem~\ref{pro_om_nuclear}) associated with a full 
Hilbert bimodule $\al M.$ with faithful left action,
satisfying the following additional properties:
\begin{itemize}
\item[(i)] The coefficient {\it C*}-algebra $\al A.$ is nuclear and 
           separable.
\item[(ii)] The left $\al A.$-action is given in terms of compact
            module operators and is non-degenerate.
\item[(iii)] $\al M.$ is a direct summand in the standard Hilbert module
             $\ot H._\al A.$ over $\al A.$.
\end{itemize}

Nuclearity of Cuntz-Pimsner algebras has been discussed recently
in the concrete Toeplitz algebra setting
(cf.~\cite{Kumjian04,pGermain02}). The paper by Kumjian is in a certain 
sense complementary to ours. In fact, he considers left $\al A.$-actions
that have a trivial intersection with the compact module operators.
Our proof uses the alternative approach in \cite{Doplicher98}. 
In particular, we analyse the structure of some spectral subspaces 
$\al O._\al M.^k$, $k\in\N_0$, that are associated to
a natural circle action. An important step in the proof
is to recognize the structure of $\al O._\al M.^1$ as an imprimitivity 
$\al O._\al M.^0$-bimodule. In this way we can apply a result by 
Brown, Green and Rieffel to the corresponding stabilizations 
and consider, roughly speaking, the Cuntz-Pimsner algebra as a crossed-product
of the zero-grade spectral subspace $\al O._\al M.^0$, which 
is shown to be a nuclear {\it C*}-algebra.

\section{Basic definitions and the main theorem}

Let $\al A.$ be a {\it C*}-algebra and $\al M.$ a Hilbert $\al A.$-module.
We denote by $\al L. ( \al M. )$ the {\it C*}-algebra 
of adjointable, right $\al A.$-module operators on ${\al M.}$ and by 
$\al K. ( \al M. ) \subseteq \al L. ( \al M. )$ 
the (closed) ideal of compact operators generated by the maps
\begin{equation}
\label{def_theta}
\theta_{\psi , \psi' } \in \al L. ( \al M. )\;,\;\; 
\psi , \psi'\in{\al M.}\;,
\quad\mr with.\quad
\theta_{\psi , \psi' } (\varphi) := 
\psi \left \langle  \psi' , \varphi \right \rangle \ , \ 
\varphi \in {\al M.} \ \ ,
\end{equation}
where $\left \langle \cdot , \cdot \right \rangle$ is the 
$\al A.$-valued scalar product defined on ${\al M.}$.

We denote by ${\ot H.}_{\al A.}$ the standard 
(countably generated) Hilbert $\al A.$-module
of sequences $(A_n)_n$ such that $\sum_n A_n^*A_n$ converges
in $\al A.$ (cf.~\cite[Example~13.1.2~(c)]{bBlackadar98}). 
It is well-known that the {\it C*}-algebra $\al K. ({\ot H.}_{\al A.})$ of 
compact, right $\al A.$-module operators on ${\ot H.}_{\al A.}$ is isomorphic to 
$\al K. \otimes \al A.$, where $\al K.$ is the {\it C*}-algebra of compact 
operators over a separable Hilbert space. 
The multiplier algebra of $\al K. ({\ot H.}_{\al A.})$ is isomorphic to
$\al L. ({\ot H.}_{\al A.})$. 
We will regard the standard module 
${\ot H.}_{\al A.}$ as a Hilbert $\al A.$-bimodule
with the obvious left $\al A.$-module action.

In this paper, we will consider Hilbert $\al A.$-modules $\al M.$
which are direct summands of ${\ot H.}_{\al A.}$. This implies that
$\al M.$ is finitely or countably generated. Conversely by Kasparov 
stabilization every countably generated Hilbert $\al A.$-module is
a direct summand of ${\ot H.}_{\al A.}$. Moreover, if $\al A.$ is unital,
then also every algebraically finitely generated 
Hilbert $\al A.$-module is a direct summand of ${\ot H.}_{\al A.}$.
The left $\al A.$-action on ${\al M.}$ is given by a
*-homomorphism $\alpha\colon \al A. \rightarrow \al L. ( \al M. )$.
In the present paper, we will assume that $\alpha$ is faithful
(in the sequel, we will identify elements of $\al A.$ with their
image in $\al L.(\al M.)$) and with image contained in $\al K.(\al M.)$. 
We will also assume that $\alpha$ is non-degenerate in the sense that
\begin{equation}
\label{eq_la}
\al A. \cdot \al K. (\al M.)
:=
{\mathrm{closed \ span}}
\left\{
A T \mid
A \in \al A.
,
T \in \al K.(\al M.)
\right\}
=
\al K.(\al M.)
\ .
\end{equation}

Note that when $\al M.$ is algebraically
finitely generated, then $\al L.(\al M.) = \al K.(\al M.)$, thus 
every left $\al A.$-action is given by compact module operators.

We denote by ${\al O.}_{\al M.}$ the Cuntz-Pimsner algebra 
associated with the Hilbert bimodule $\al M.$ (cf.~\cite{Pimsner97}).
Recall that ${\al O.}_{\al M.}$ is generated as a {\it C*}-algebra by
$\al M.$ and $\al A.$ satisfying the relations
\begin{eqnarray}
\psi^* \psi' &=& \left \langle \psi , \psi' \right\rangle 
  \ , \quad \psi , \psi' \in \al M.
\label{def_pimsner}\\
A \psi     &:=&  \alpha (A) \psi \ ,
 \quad A \in \al A.  \;.\label{def_pimsner1}
\end{eqnarray}
Note that (\ref{def_pimsner}) implies 
$\psi' \psi^* \varphi = \theta_{\psi' , \psi} (\varphi)$, 
$\psi' , \psi^* , \varphi \in {\al M.}$. 
Therefore one has the natural identification
\begin{equation}
\label{eq_compact_operators}
\theta_{\psi' , \psi} = \psi' \psi^* \ \ .
\end{equation}

Moreover, if
${\al M.}^r := {\al M.} \otimes_{\al A.} \cdots \otimes_{\al A.} {\al M.}$,
$r\in\N$, 
is the $r$-fold tensor product with coefficients in $\al A.$,
then there is an identification
\begin{equation}
\label{eq_tens}
{\al M.}^r \simeq {\mathrm{closed \ span}}\left\{ 
\psi_1\cdot \ldots \cdot\psi_r \in {\al O.}_{\al M.} \mid
\psi_k \in {\al M.} , k = 1 , \ldots , r \right\} \,.
\end{equation}

There is a natural action of the circle
$\T:=\{z\in\C\mid |z|=1\}$ on the Cuntz-Pimsner algebra
given by
\begin{equation}
\label{def_circle_action}
\delta\colon \T \rightarrow \aut {\al O.}_{\al M.} \ \ , \ \ 
\delta_z (\psi) := z \psi \ \ , \ \ 
z \in \T , \ \psi \in {\al M.} \,.
\end{equation}
We denote by 
\begin{equation}\label{SpecS1}
{\al O.}^k_{\al M.} := \left\{ 
T \in {\al O.}_{\al M.} \mid \delta_z (T) = z^k T \ , \ z \in \T
\right\}
\ , \ 
k \in \Z \ ,
\end{equation}
the spectral subspaces associated to the circle action. 
In particular, ${\al O.}^0_{\al M.}$ is the closed span of 
elements of the form $\psi' \psi^*$, $\psi , \psi' \in \al M.^r$,
$r \in \N$. From Eq.~(\ref{eq_la}) we have that 
\[ 
 \al A. \cdot \al M. :=\mathrm{closed~span} \left\{ A \psi \mid A \in \al A. , \psi \in \al M. \right\}
                     =\al M.\,.
\]
(Use the fact that every $\psi \in \al M.$ is of
the form $\psi = T \psi'$ for some $T \in \al K.(\al M.)$,
$\psi' \in \al M.$, cf.~\cite[Lemma 1.3]{Blanchard96}).
Therefore $\al A. \cdot \al M.^r=\al M.^r$, $r\in\N$, and
\begin{equation}
\label{eq_la2}
\al A. \cdot \al O.^0_{\al M.}
=
{\mathrm{closed~span}}
\left\{
A T \mid
A \in \al A. \ , \
T \in \al O.^0_{\al M.}
\right\}
=
\al O.^0_{\al M.} \ .
\end{equation}

In order to discuss 
universality properties of the Cuntz-Pimsner algebra we give the 
following definitions (cf.~\cite[\S 2]{Doplicher98})

\begin{defi} \label{def_hba} 
Let $\al A. \subset \al B.$ be an inclusion of {\it C*}-algebras. A {\bf Hilbert 
$\al A.$-bimodule in} $\al B.$ is a closed vector space 
${\al M.} \subset \al B.$ satisfying
\begin{itemize}
\item[(i)] $A \psi \in {\al M.}$, $\psi A \in {\al M.}$ and
$\psi^* \psi' \in \al A.$ for every $A \in \al A.$, $\psi , \psi' \in {\al M.}$.
\item[(ii)] For any $A\in\al A.$ with $A \psi = 0$, $\psi \in \al M.$, one has $A = 0$.
\end{itemize}
We say that $\al M.$ is full if for every $A \in \al A.$ there are 
$\psi ,\psi' \in {\al M.}$ such that $A = \psi^* \psi'$.
\end{defi}

Now, there is a natural identification
\[
\al M. \al M.^* :=
{\mathrm{closed \ span}}
\left\{
\psi' \cdot\psi^* \mid \psi , \psi' \in \al M. \subset \al B.
\right\}
\simeq
\al K. (\al M.) \ ,
\]
hence $\al M. \al M.^* $ may be regarded as a {\it C*}-subalgebra of $\al B.$. 
We say that $\al M.$ has support $\1$ if there exists a sequence
$\left\{ \psi_n \right\} \subset$ $\al M.$ such that
$\left\{  U_\mt N. := \sum_{n=1}^\mt N. \psi_n \psi_n^* \right\}_\mt N.$
is an approximate unit for $\al B.$ (recall that by assumption the 
{\it C*}-algebras are separable.)

\begin{rem}
Note that if $\al M. \subset$ $\al B.$ has support $\1$ and if there 
are elements $A \in \al A.$, $T  \in\al M. \al M.^*$ satisfying
$A \psi = T \psi$ for every $\psi \in \al M.$, then 
\begin{equation}
\label{tea}
T = A \,.
\end{equation}
In fact, since $T \psi_n = A \psi_n$ for every $n \in \N$, we have
$T U_\mt N. = A U_\mt N.$ for every $N \in \N$ and therefore 
\[
 T =\lim_\mt N. T U_\mt N. =\lim_\mt N. A U_\mt N. = A\,.
\]
\end{rem}

The following result is just a translation of \cite[Theorem~3.12]{Pimsner97}
in terms of Hilbert bimodules in {\it C*}-algebras. We note 
explicitely that Eq.~(\ref{tea})
is equivalent to condition (4) in the above-cited theorem.

\begin{pro} \label{prop_uni}
Let $\al A. \subset \al B.$ be an inclusion of unital {\it C*}-algebras, 
${\al M.} \subset \al B.$ a full
Hilbert $\al A.$-bimodule in $\al B.$ 
with support $\1$. Then there is a canonical morphism 
${\al O.}_{\al M.} \rightarrow \al B.$.
\end{pro}

\begin{rem} \label{rem_partial_isometry}
Examples of the above universality property can be found in the Cuntz-Pimsner 
algebra itself. 
In fact, if $\al M.$ is algebraically finitely generated
and $\al A.$ is unital, then it follows from the
definition of the Cuntz-Pimsner algebra that $\sum_n \psi_n \psi_n^* =$
$\1$ for every finite set $\left\{ \psi_n \right\}$ of (normalized) generators 
of $\al M.$. If $\al M.$ is countably generated as a right $\al A.$-module,
then there are elements $\left\{ \psi_n \right\}_{n=1}^\infty\subset\al M.$
such that $U_\mt N. :=\sum_{n=1}^\mt N. \psi_n \psi_n^*$
is an approximate unit for $\al K.(\al M.)$, hence also for 
$\al A.$ which may be regarded as a {\it C*}-subalgebra of $\al K.(\al M.)$
(see for example \cite{Brown77a} or p.~266 in \cite{Doplicher98}). 
Finally, since $\al O._\al M.$ is generated as a {\it C*}-algebra by $\al M.$ and $\al A.$, 
we conclude that $\left\{ U_\mt N. \right\}_\mt N.$ is an approximate unit for
${\al O.}_{\al M.}$, so that $\al M.$ has support $\1$ in $\al O._\al M.$.
\end{rem}

\begin{rem} \label{rem_ok}
The {\it C*}-algebra $\al K.$ of compact 
operators over a separable Hilbert space
is clearly a Hilbert $\al K.$-bimodule
with left and right actions defined by multiplication
and scalar product 
$
\left \langle V , V'  \right \rangle := V^* V'
$, $V , V' \in \al K.$.
Let $\al M.$ be a Hilbert $\al A.$-bimodule.
We consider the external tensor product of Hilbert bimodules
$\al M. \widehat \otimes \al K.$, which has
a natural structure of Hilbert $(\al A. \otimes \al K.)$-bimodule,
and denote by $\al O._{\al M. \widehat \otimes \al K.}$
the associated Cuntz-Pimsner algebra.
Then, there is a natural identification
$
\al O._{\al M. \widehat \otimes \al K.} 
\simeq
\al O._\al M. \otimes \al K.
$
defined by the map 
$
\psi \widehat \otimes V
\mapsto 
\psi \otimes V
$,
$\psi \in \al M.$, $V \in \al K.$.
\end{rem}

\begin{pro}\label{rem_om0_nuc}
Let $\al M.$ be a Hilbert $\al A.$-bimodule satisfying properties 
(i)-(iii) of the introduction and let $\al O._\al M.$ be the corresponding
Cuntz-Pimsner algebra. Then the zero grade {\it C*}-algebra
$\al O._\al M.^0$ (cf.~(\ref{SpecS1})) is nuclear.
\end{pro}
\begin{beweis}
First we prove that every $\al K. (\al M.^r)$, $r \in \N$, is a nuclear 
{\it C*}-algebra. Since $\al M.$ is countably generated (or algebraically 
finitely generated if $\al A.$ is unital), we conclude that $\al M.^r$ 
is countably generated (or algebraically finitely generated if $\al A.$ is
unital). By Kasparov stabilization, we obtain that $\al M.^r$ is a direct
summand of $\ot H._\al A.$. This implies that $\al K. (\al M.^r)$
is a corner of the nuclear {\it C*}-algebra 
$\al K.(\ot H._\al A.) \simeq$ $\al A. \otimes \al K.$ and, therefore,
$\al K.(\al M.^r)$ is nuclear.
Now, for every $r \in \N$ there is an embedding
\[
i_r \colon \al K. (\al M.^r) \rightarrow \al L. (\al M.^{r+1})
\ \ , \ \
i_r (T) \ \psi \varphi := (T \psi )  \varphi
\ ,
\]
$\psi \in \al M.^r$, $\varphi \in \al M.$
(in the usual tensor notation the embedding is given simply
by $i_r(T) :=T \otimes \1$,
where $\1$ is the identity of $\al L.(\al M.)$). 
Next we show that
if the image of the left $\al A.$-action is contained in $\al K.(\al M.)$,
then $i_r(T) \in\al K. (\al M.^{r+1})$. First put
$T := \psi_1 \psi_2^*$, $\psi_1$, $\psi_2 \in$ $\al M.^r$
so that, $i_r(T) \psi \varphi =\psi_1 \psi_2^* \psi \varphi$. 
By \cite[Lemma 1.3]{Blanchard96},
there is a decomposition $\psi_1 = \psi_0 A_0$ for some 
$A_0 \in \al A.$, $\psi_0 \in \al M.^r$. 
Moreover, since the left $\al A.$-module
action is given by compact operators,
there exist $\varphi_0,\varphi_1 \in\al M.$ with
$A_0 = \varphi_0 \varphi_1^*$, hence 
$i_r(T) =$ $\psi_0 \varphi_0 \varphi_1^* \psi_2^*$ is an
element of $\al K. (\al M.^{r+1})$.
Therefore the zero grade algebra $\al O._{\al M.}^0$ is an inductive 
limit 
\[
\al O._{\al M.}^0 
\ = \ 
\lim_{\longrightarrow}           
( \al K. (\al M.^r)  ,  i_r )
\ .
\]
Since every $\al K.(\al M.^r)$ is nuclear, we conclude that 
$\al O._{\al M.}^0$ is nuclear.
\end{beweis}

We can now prove our main theorem.

\begin{teo}
\label{pro_om_nuclear}
Let $\al A.$ be a nuclear and separable {\it C*}-algebra. Assume that
${\al M.}$ is a full and non-degenerate
Hilbert $\al A.$-bimodule with faithful left $\al A.$-action and satisfying:
\begin{itemize}
\item[(i)] $\al M.$ is a direct summand (as a right Hilbert $\al A.$-module) of the 
           standard Hilbert module ${\ot H.}_{\al A.}$ over $\al A.$.
\item[(ii)] The left $\al A.$-action on $\al M.$ is given in terms of compact 
            module operators.
\end{itemize}
Then the Cuntz-Pimsner algebra ${\al O.}_{\al M.}$ is nuclear.
\end{teo}
\begin{beweis}
Consider the spectral subspace $\al O.^1_{\al M.} \subset {\al O.}_{\al M.}$
which was introduced in (\ref{SpecS1}) and note that it has a natural structure
as a $\al O.^0_{\al M.}$-bimodule. In fact, take left and right multiplication
by elements of $\al O.^0_{\al M.} \subset {\al O.}_{\al M.}$ and define the
$\al O.^0_{\al M.}$-valued scalar product by
\[
\left \langle T , T' \right \rangle := T^* T' 
\ \ , \ \
T,T' \in \al O.^1_{\al M.} \ .
\]
We show next that $\al O.^1_{\al M.}$ is a full Hilbert 
$\al O.^0_{\al M.}$-module: by Eq.~(\ref{eq_la2}) we have that
$\al A. \cdot \al O.^0_{\al M.} = \al O.^0_{\al M.}$, so that
any $T \in\al O.^0_{\al M.}$ can be written as
$T = AT'$ for some $A \in \al A.$ and $T' \in \al O.^0_{\al M.}$.
Since $\al M.$ is full, there are $\psi, \psi' \in \al M.$
such that $A = \psi^* \psi'$ and therefore
\[
T = AT'=\psi^* \psi' T' =\left \langle \psi , \psi' T ' \right \rangle
       \in\langle\al O.^1_{\al M.},\al O.^1_{\al M.}\rangle\,.
\]
Therefore $\al O.^1_{\al M.}$ is full as a Hilbert $\al O.^0_{\al M.}$-module. 

Denote by 
\[
\ot O.:=\al O._{\al O.^1_{\al M.}}
\]
the Cuntz-Pimsner algebra associated with the $\al O.^0_{\al M.}$-bimodule
$\al O.^1_{\al M.}$. It is enough to show that $\ot O.$ is nuclear since
$\ot O.$ is isomorphic the original Cuntz-Pimsner algebra ${\al O.}_{\al M.}$.
In fact, ${\al O.}_{\al M.}^1$ is a 
Hilbert ${\al O.}_{\al M.}^0$-bimodule in ${\al O.}_{\al M.}$ with 
support $\1$ so that by Proposition~\ref{prop_uni} there exists a monomorphism 
$I\colon \ot O. \hookrightarrow {\al O.}_{\al M.}$. Moreover, since 
$\al M. \subset {\al O.}_{\al M.}^1$, it follows that $I$ is surjective
(cf.~\cite[Theorem~2.5]{Pimsner97}).

To show the nuclearity of $\ot O.$ we need to exploit 
the additional structure of the bimodule $\al O.^1_{\al M.}$.
Note first that $\al O.^1_{\al M.}$ also carries a natural 
left $\al O.^0_{\al M.}$-valued scalar product given by
\[
\left \langle T' , T \right \rangle_l := T' T^* 
\ \ , \ \
T,T' \in \al O.^1_{\al M.} \ .
\]
By a similar argument as before we have that 
$\al O.^1_{\al M.}$ is also full w.r.t.~$\langle\cdot,\cdot\rangle_l$,
i.e.~$\al O.^0_{\al M.} = \left
  \langle \al O.^1_{\al M.} , \al O.^1_{\al M.} \right \rangle_l$.
We conclude that $\al O.^1_{\al M.}$ can also be interpreted as an 
$\al O.^0_{\al M.}$-imprimitivity bimodule
(cf.~\cite[Section~3.1]{bRaeburn98}). 
 
Put
\[
 \al J.:=\al O.^1_{\al M.} \widehat \otimes \al K.
 \quad\mr as~well~as.\quad 
 \al B.:= \al O.^0_{\al M.} \otimes \al K.\,,
\]
where $\widehat\otimes$ denotes the external tensor product of 
Hilbert bimodules and $\al K.$ is the {\it C*}-algebra of compact 
operators over a separable Hilbert space.
Then, $\al J.$ is an imprimitivity $\al B.$-bimodule
(cf.~Remark~\ref{rem_ok}). Since $\al B.$ is a stable and separable
{\it C*}-algebra we obtain from Corollary~3.5 in \cite{Brown77}
that there is an isomorphism of Hilbert bimodules
\begin{equation}\label{beta0}
 \beta_0\colon \al J.\to \al B.\,,
\end{equation}
where $\al B.$ is considered as a bimodule over itself with 
multiplication as right action and the left action
being specified by a suitable automorphism 
$\theta\in\mr Aut.\al B.$. 
The isomorphism $\beta_0$ extends to 
an isomorphism of the corresponding Cuntz-Pimsner algebras.
Moreover, $\beta_0$ also extends to an isomorphism of the
associated multiplier algebras.
Hence we have
\begin{equation}\label{IsomCross}
\beta_0 \colon \al O._\al J.\to \al O._{\al B.}\simeq \al B. \rtimes_\theta \Z\,,
\end{equation}
where for the last isomorphism with the crossed product we
use the results in \cite[Chapter~1]{Pimsner97}.

Let $E\in\al K.$ be a minimal projection. Then using Eq.~(\ref{beta0})
we may define a monomorphism
\begin{equation}
  \label{beta}
  \beta\colon\al O.^1_{\al M.}\to \al B.\quad \mr by~means~of.\quad
  T\mapsto \beta_0(T\otimes E)\,.
\end{equation}
Note that the image of $\al O.^1_{\al M.}$ generates
$\al B.$ as a $\al B.$-bimodule and that 
(by universality) $\beta$ can be extended to a monomorphism
\begin{equation}
  \label{Mono}
  \beta\colon\ot O.\to \al O._{\al B.} \,,
\end{equation}
where $\ot O.:=\al O._{\al O.^1_{\al M.}}$ was introduced in the beginning
of the proof. 
Since $\al O.^0_{\al M.}$ is nuclear (cf.~Proposition~\ref{rem_om0_nuc}) 
we have that $\al B.$ is nuclear and the same is true for the crossed product
$\al B. \rtimes_\theta \Z$ (\cite[Theorem~15.8.2]{bBlackadar98}).
From Eq.~(\ref{IsomCross}) we obtain that $\al O._\al J.$ and $\al O._{\al B.}$ 
are nuclear {\it C*}-algebras.

Finally, we turn our attention to the Cuntz-Pimsner algebra $\ot O.$.
We will conclude the proof by showing that this algebra is a corner of the nuclear
algebra $\al O._{\al B.}$: By Remark~\ref{rem_ok} we may 
identify $\al O._\al J.\simeq\ot O.\otimes\al K.$ and using
(\ref{IsomCross}) we conclude that $\{\beta_0(T\otimes V)\mid
T\in\ot O.,V\in\al K.\}$ is total in $\al O._{\al B.}$. 
Let us now consider the identity $\1$ of the multiplier algebra $M(\ot O.)$.
Then, $\1 \otimes E \in$ $M (\ot O. \otimes \al K.)$ and we define the
projection 
$E_\beta := \beta_0(\1\otimes E) \in$ $M(\al O._{\al B.})$.
For $T \in \ot O.$ we have that
$\beta (T) = E_\beta \beta (T) E_\beta$.
On the contrary, take $B = \sum_{i=1}^n \beta_0(T_i \otimes V_i)\in$
$\al O._\al B.$, $T_i \in \ot O.$, 
$V_i\in\al K.$, $i=1,\ldots,n$. Then
\[
 E_\beta B E_\beta=\sum_i E_\beta\beta_0(T_i \otimes V_i)E_\beta
                  =\sum_i \beta_0 ( T_i \otimes EV_iE ) 
                  =\beta_0 \left( \left(\sum_i z_i T_i\right) \otimes E \right)\,,
\]
where $z_i\in\C$ are given by $EV_iE=z_iE$ (recall that $E$ is minimal). 
We conclude that $\ot O.$ is isomorphic to the corner
$E_\beta \al O._\al B. E_\beta$. Thus, $\ot O.$ is nuclear.
\end{beweis}


\section{Outlook}

Cuntz-Pimsner algebras provide an important family of examples in the theory 
operator algebras. Moreover, these algebras, which generalise Cuntz algebras, 
appear naturally in the extension of Doplicher-Roberts superselection
theory (cf.~\cite{Doplicher90})
to the case where the observable algebra has a nontrivial center $\al Z.$
(see e.g.~\cite{Baumgaertel05,Lledo04a,Lledo97b}).
In this context the category of canonical endomorphisms
is isomorphic to a category whose objects are free $\al Z.$-bimodules.
The Cuntz-Pimsner algebras associated to these $\al Z.$-bimodules
generate a {\it C*}-algebra $\al F.$ on which one can realise concretely
the dual of a compact group
(cf.~\cite{pZioNando06b} and references cited therein). 
The class of algebras considered in the present paper contain the Cuntz-Pimsner 
algebras that appear in this application. In certain
special cases, tensor products of these algebras may also appear in
concrete models.

\paragraph{Acknowledgments}
We are grateful to the DFG-Graduiertenkolleg 
``Hierarchie und Symmetrie in mathematischen Modellen''
for supporting a visit of E.V.~to the RWTH-Aachen University.
E.V. was also partially supported by the European Network 
``Quantum Spaces - Noncommutative Geometry" HPRN-CT-2002-00280.


\end{document}